\documentclass[a4paper,11pt]{amsart}
\usepackage[cp1251]{inputenc}
\usepackage[english]{babel}
\usepackage{amsmath,amssymb,euscript,amsthm,amsfonts}
\newtheorem{lem}{Lemma}
\newtheorem{thm}{Theorem}
\theoremstyle{definition}
\newtheorem{defn}{Definition}
\theoremstyle{remark}
\newtheorem{remk}{Remark}
\newtheorem{expl}{Example}

\begin{document}

\renewcommand{\proofname}{Proof}
\makeatletter \headsep 10 mm \footskip 10 mm
\renewcommand{\@evenhead}%
{\vbox{\hbox to\textwidth{\strut \centerline{{\it Nikolaj
Glazunov}}} \hrule}}

\renewcommand{\@oddhead}%
{\vbox{\hbox to\textwidth{\strut \centerline{{\it
Critical lattices, elliptic curves and their possible dynamics
%If the title is too long for a line, please, write our its shorten
% variant
}}} \hrule}}

%Papers for the Proceedings of the Third Vorono\"\i \ Conference of
%the Number Theory and Spatial Tessellations

\begin{center}
{\Large\bf CRITICAL LATTICES, ELLIPTIC CURVES AND THEIR POSSIBLE DYNAMICS}
\end{center}

\begin{center}
NIKOLAJ  GLAZUNOV\ (Ukraine)
\end{center}

\vskip 15pt

\begin{quote}
 We present a combinatorial geometry and dynamical systems framework
 for the investigation and proof of the Minkowski conjecture about critical determinant
 of the region $ |x|^p + |y|^p < 1, \; p > 1. $
 The application of the framework may drastically reduce the
 investigation of sufficiently smooth real functions of many
 variables.
 Incidentally, we establish connections between critical lattices,
 dynamical systems  and elliptic curves.
\end{quote}

\vskip 20pt

\section{Introduction}
 Vorono{\"i} [\ref{V:CW}] have showed that a lattice is extreme if
 and only if it is both perfect and eutactic.  The notions of the
 critical lattice are partial case of extreme lattice.
 Let $\mathcal D$  be a set in $n-$dimensional real space ${\bf R}^n.$ A
lattice $\Lambda $
 is the {\it admissible} for the set ${\mathcal D}$
($ {\mathcal D}-{\mathit  admissible}$) if
${\mathcal D} \bigcap \Lambda = \emptyset $ or $0.$
The infimum $\Delta(\mathcal D) $ of determinants of all lattices
admissible for ${\mathcal D}$ is called {\it the critical determinant}
of $\mathcal D. $
  A lattice $\Lambda $ is the {\it critical}
if $ d(\Lambda) = \Delta({\mathcal D})$ [\ref{C:GN}].
 Let now $ D_p \subset {\bf R}^2 = (x,y), \ p > 1 $ be the 2-dimensional
region: $ |x|^p + |y|^p < 1 . $ Let $\Delta(D_p) $ be the critical
determinant of the region. H. Minkowski [\ref{Mi:DA}] have raised
a question about critical determinants and critical lattices of
regions $ D_p$ under varying $p > 1 $. Let $\Lambda_{p}^{(0)} $
and $\Lambda_{p}^{(1)} $ be two ${\mathcal D}-$admissible lattices each
of which contains three pairs of points on the boundary of $ D_p$
and such that $ (0,1) \in \Lambda_{p}^{(0)}, \;
(-2^{-1/p},2^{-1/p}) \in \Lambda_{p}^{(1)},  $ (under these
conditions the lattices are unique defined). Investigations by H.
Minkowski, L. Mordell, C. Davis, H. Cohn, A. Malishev and another
researches (please, see the bibliography at [\ref{GGM:PM}]) in
terms of lattices gave [\ref{GGM:PM}]:
\begin{thm}
 {\it $$\Delta(D_p) = \left\{
                   \begin{array}{lc}
    d(\Lambda_{p}^{(1)}), \; 1 < p \le 2, \; p \ge p_{0},\\
    d(\Lambda_{p}^{(0)}), \;  2 \le p \le p_{0};\\
                     \end{array}
                       \right.
                           $$
here $p_{0}$ is a real number that is defined unique by conditions
$d(\Lambda_{p_{0}}^{(0)}) = d(\Lambda_{p_{0}}^{(1)})$, $2,57 \le
p_{0} \le 2,58. $}
\end{thm}

The Minkowski conjecture about critical determinant of the region
$ D_p $ can be formulated as the problem of minimization on moduli
space ${\mathcal M}$ of admissible lattices of the region $ D_p $
(see below). This moduli space is a differentiable manifold. We
will call it the {\it Minkowski's moduli space}. The {\it tangent
bundle} of a differentiable manifold $ M$, denoted by $ TM $, is
the union of the tangent spaces at all the points of $ M .$ Recall
that a {\it vector field} on a smooth manifold $M$ is a map $F: M
\rightarrow TM $  which satisfies $ p \circ F = id_{M} $, where
$p$ is the natural projection  $ TM \rightarrow M. $ By its
definition a vector field is a {\it section} of the bundle $ TM .$
There are many classes of dynamical systems on the Minkowski's
moduli space. Each vector field defines a dynamical system on $M$.
But we do not consider these dynamical systems in the paper. There
are bundles of sufficiently smooth real functions  on ${\mathcal
M} $ and sections of the bundles define discrete dynamical systems
on ${\mathcal M}$.  Below we will define interval extension of
${\mathcal M} $ and  interval sheaves on the interval extension
${\mathcal IM}$ of ${\mathcal M}$ and  interval dynamical systems
on ${\mathcal IM}.$ Still one class of dynamical systems appear
from lattices. After complexification of ${\bf R}^2$ and lattices
$\Lambda_{p}^{(0)}$ and $\Lambda_{p}^{(1)}$ these lattices define
two classes of elliptic curves. C. Deninger [\ref{D:AN}]  and authors of
the paper [\ref{DEMW:DS}] have discussed possible dynamics of elliptic
curves. At the last section of the paper we will investigate
elliptic curves of above-mentioned classes and their possible
dynamics. In order to make the paper easy to read we will
illustrate definitions of the basic notions by examples.

%\section{Minkowski conjecture about critical determinant of
 %$|x|^p + |y|^p <1, p > 1,$  domain}

%\begin{defn}
%Statement of definition
%\end{defn}

%\begin{lem}
%Statement of Lemma
%\end{lem}

%\begin{proof}
%Lemma's proof
%\end{proof}

%\begin{thm}
%Statement of Theorem
%\end{thm}

%\begin{proof}
%Theorem's proof
%\end{proof}

%\begin{remk}
%Statement of remark
%\end{remk}

%\subsection{Type Changes and {\subsecit Special} Characters}

%\section{History of the investigation of the Minkowski conjecture}

\section{Analytical formulation of the Minkowski's conjecture}
Recall the analytic formulation of Minkowski's conjecture
[\ref{Mi:DA}], [\ref{D:NC}], [\ref{Co:MC}]. Below we use notations
from [\ref{GGM:PM}], [\ref{G:RM}].
 Let
$$ \Delta(p,\sigma) = (\tau + \sigma)(1 + \tau^{p})^{-\frac{1}{p}}
  (1 + \sigma^p)^{-\frac{1}{p}},  $$
be the function defined in the domain
 $$ D_{p}: \; \infty > p > 1, \; 1 \leq \sigma \leq \sigma_{p} =
 (2^p - 1)^{\frac{1}{p}}, $$
of the $ \{p,\sigma\} $ plane, where $\sigma$ is some real
parameter; $\;$ here $ \tau = \tau(p,\sigma) $ is the function
uniquely determined by the conditions $$ A^{p} + B^{p} = 1, \; 0
\leq \tau \leq \tau_{p}, $$ where $$ A = A(p,\sigma) = (1 +
\tau^{p})^{-\frac{1}{p}} - (1 + \sigma^p)^{-\frac{1}{p}}, $$ $$ B
= B(p,\sigma) = \tau(1 + \sigma^p)^{-\frac{1}{p}}       + \sigma(1
+ \tau^{p})^{-\frac{1}{p}}, $$ $\tau_{p}$ is defined by the
equation $ 2(1 - \tau_{p})^{p} = 1 + \tau_{p}^{p}, \; 0 \leq
\tau_{p} \leq 1. $
\begin{expl}
  Critical lattices for $p = 2.$ In the case there are two
  critical lattices: $\Lambda_{2}^{(0)}$ and $\Lambda_{2}^{(1)}.$
  The lattice $\Lambda_{2}^{(0)}$ has basis $\omega_1 = (1, 0), \;
  \omega_2 = (1/2, {\sqrt 3}/2).$ \\
   The lattice $\Lambda_{2}^{(1)}$ has basis $\omega_1 = (-2^{-1/2}, 2^{-1/2}), \;
  \omega_2 = (\frac{{\sqrt{2 - {\sqrt 3}}}}{2}, \frac{{\sqrt{2 + {\sqrt 3}}}}{2}).$
  \end{expl}
  \begin{expl}
  More generally, for $2 \le p \le p_{0}, $ the critical lattice
  $\Lambda_{p}^{(0)}$ has the basis $\omega_1 = (1, 0), \;
  \omega_2 = (1/2, \sigma_{p}/2).$
  \end{expl}
\begin{defn}
 In the notations the surface
$$ \Delta - (\tau + \sigma)(1 + \tau^{p})^{-\frac{1}{p}}
  (1 + \sigma^p)^{-\frac{1}{p}} = 0,  $$
  in $3-$dimensional real space with coordinates $(\sigma,p,\Delta)$ is
called the   {\it Minkowski's moduli space}.
\end{defn}
 {\bf Minkowski's analytic conjecture (MA):}
\quad
 For any real $p$ and $\sigma$ with conditions
$ p > 1, \ p \ne 2, \ 1 < \sigma < \sigma_{p} $ \\ $$
\Delta(p,\sigma) > \Delta(D_p) = \min
(\Delta(p,1),\Delta(p,\sigma_p)). $$ \\

\section{Interval Cellular Covering}

\begin{defn}
  For any $n$ and any $j, \; 0 \leq j \leq n$, an $j-${\it dimensional
interval cell}, or $j-${\it I-cell}, in ${\bf R}^n$ is a subset
$Ic$ of ${\bf R}^n$ such that (possibly, after permutation of
variables) it has the form   \\

$ Ic = \{x \in {\bf R}^n: {\underline a}_{i},{\overline a}_{i},
r_{k} \in {\bf R} :    \\
 \; {\underline a}_{i} \le x_{i} \le {\overline a}_{i}, 1 \le i \le j, \\
 \; x_{j + 1} = r_{1}, \cdots, x_{n} = r_{n-j} \} \; . $
Here $  {\underline a}_{i} \le {\overline a}_{i}. $
\end{defn}
If $j = n$ then we have an $n-$dimensional interval vector. Let
${\mathcal P}$ be the hyperplane that contains $Ic.$ These is the well
known fact:
\begin{lem}
 The dimension of $Ic$ is equal to the minimal dimension of hyperplanes
that contain $Ic.$
\end{lem}

 Let ${\mathcal P}$ be the such hyperplane, $Int \; Ic$ the set of interior
points of $Ic$ in ${\mathcal P}, \; Bd \; Ic = Ic \setminus Int \;
Ic.$ For $m-$dimensional I-cell $Ic$ let $d_{i}$ be an
$(m-1)-$dimensional I-cell from $Bd \; Ic.$ Then $d_{i}$ is called
an $(m-1)-$dimensional {\it face} of the I-cell $Ic.$

\begin{defn}
 Let $D$ be a bounded set in ${\bf R}^n.$ By interval cellular covering
$Cov$ we will understand any finite set of $n-$dimensional I-cells
such that their union contains $D$ and adjacent I-cells are
intersected by their faces only. By $ {\mid Cov \mid}$ we will
denote the union of all I-cells from $Cov.$
\end{defn}

  Let $Cov$ be the interval covering. By its  {\it subdivision}
we will understand an interval covering $Cov^{`}$ such that $\mid
Cov \mid = \mid Cov^{`} \mid$ and each I-cell from $Cov^{`}$ is
contained in an I-cell from $Cov.$ In the paper we will consider
mainly bounded horizontal and vertical strips in ${\bf R}^2,$
their interval coverings and subdivisions.

\section{Some Categories and Functors of Interval Mathematics}

 Let $ {\bf X} =
({\bf x}_{1}, \cdots,{\bf x}_{n}) = ([{\underline x}_{1},
{\overline x}_{1}], \cdots, [{\underline x}_{n}, {\overline
x}_{n}] $ be the $n-$dimensional real interval vector with $
{\underline x}_{i} \leq x_{i} \leq {\overline x}_{i} $
 ("rectangle" or "box").
Let $f$ be a real continuous function  of  $n$ variables that is
defined on ${\bf X}$. The {\it interval evaluation} of $f $ on the
interval ${\bf X} $ is the interval $[{\underline f}, {\overline
f}] $ such that for any $ x \in {\bf X}, \; f(x) \in [{\underline
f}, {\overline f}]. $ The interval evaluation is called {\it
optimal} [\ref{AH:IC}] if $ {\underline f} = \min f,$ and $
{\overline f} = \max f $ on the interval {\bf X}. Let $ Of-$ be
the optimal interval evaluation of $f$  on ${\bf X}$.

\begin{defn}
 The pair $({\bf X}, Of)$ is
called the {\it interval functional element.} If $Ef$ is an
interval that contains $Of$ then we will call the pair $({\bf
X},Ef)$ {\it the extension} of $({\bf X}, Of)$ or  $eif-${\it
element.}      \\ Let $f$ be the constant signs function on ${\bf
X}.$ If $f > 0$ (respectively $f < 0$) on ${\bf X}$ and $Of > 0$
(respectively $Of < 0$) then we will call $({\bf X}, Of)$ {\it the
correct interval functional element} (shortly $c-${\it element}).
\end{defn}

%\begin{expl}
% An interval evaluation and the optimal interval evaluation.
%\end{expl}

More generally we will call {\it the correct interval functional
element} an extension $({\bf X},Ef)$ of $({\bf X}, Of)$ that has
the same sign as $Of.$

A set of intervals with inclusion relation forms a category
${\mathcal CIP}$ of preorder [\ref{G:IA}].
\begin{defn}
  A contravariant functor from ${\mathcal CIP}$ to the category of sets
is called the interval presheaf.
\end{defn}

%\begin{expl}
% Objects of ${\mathcal CIP}$ and action of the interval presheaf
% of objects and morphisms.
%\end{expl}

  For a finite set $FS = \{{\bf X}_{i}\}$ of $m-$dimensional intervals in
${\bf R}^n, \; m \le n,$ the union $V$ of the intervals forms a
piecewise-linear manifold in ${\bf R}^n.$ Let $G$ be the graph of
the adjacency relation of intervals from $FS$. The manifold $V$ is
connected if $G$ is a connected graph. In the paper we are
considering connected manifolds. Let $f$ be a constant signs
function on ${\bf X} \in FS.$ The set $\{({\bf X}_{j},Of)\}$ of
$c-$elements (if exists) is called {\it a constant signs
continuation of} $f$ on $\{{\bf X}_{j}\}.$ If $\{{\bf X}_{j}\}$ is
the maximal subset of $FS$ relatively a constant signs function
$f$ then $\{({\bf X}_{j},Of)\}$ is called the {\it constant signs
continuation of} $f$ on $FS$.

\section{Interval iterative processes}

In the section a "dynamical system" is a continuous map $ T: \; X
\rightarrow X$ (or continuous flow or semiflow $ \phi_{t}: \; X
\rightarrow X$) on a compact metric space. The recent
investigations of iterated polynomial maps which can be considered
as dynamical systems (see, for instance [\ref{Mil:IC}] and
references in the paper) has provoked an interest in the range of
applicability of these methods for transcendental iterated maps
and for their interval extensions. Interval mathematics offers a
rigorous approach to computer investigation of mathematical
models. Interval methods is a kind of numerical methods with
automatic result verification. Under the investigation and proof
of the Minkowski conjecture we have to compute expressions $
\Delta_{\sigma}^{'} \; , \; \Delta_{\sigma^{2}}^{''} \; , \;
\Delta_{p}^{'} \; , \;
 \Delta_{\sigma p}^{''} \; , \; \Delta_{\sigma^{2}p}^{'''} \; $
 and their interval extensions. These expressions are represented
 in terms of a sum of derivatives of "atoms" $ s_{i} = \sigma^{p-i},
\; t_{i} = \tau^{p-i}, \; a_{i} = (1 +
\sigma^{p})^{-i-\frac{1}{p}}, \; b_{i} = (1 +
\tau^{p})^{-i-\frac{1}{p}}, \; A = b_{0} - a_{0}, \; B = \tau
b_{0} + \sigma a_{0}, \; \alpha_{i} = A^{p-i}, \; \beta_{i} =
B^{p-i} \; ( i = 0, 1, 2, \ldots).$ Let $ D $ be a subdomain of  $
D_{p}. $
%Under evaluation in $ D $  mentioned functions
The domain is covered by rectangles of the form $$ {\bf X} =
[{\underline p}, {\overline p}; \; {\underline \sigma}, {\overline
\sigma}].$$ Let $f$ be one of mentioned functions. The
$eif-$element $({\bf X},Ef)$ is represented
 in terms of $
{\underline p}, \; {\overline p}, \; {\underline \sigma}, \;
{\overline \sigma}, \; {\underline \tau}, \; {\overline \tau}, \;
; \; $ here the bounds $ \; {\underline \tau},
 \; {\overline \tau}, \; $ are obtained with the help of some interval iteration
processes. In this section we give formulas for one of the
interval iterative processe on Minkowski's moduli space.
Computations of these iterative processes and their interval
extensions in various floating points and intervals were produced
%[\ref{Gl90}], [\ref{Gl89}].       \\
 Let $ f: {\bf R}^{2}
\rightarrow {\bf R} $ be a transcendental map which is a
superposition of rational functions, exponential functions and
logarithms. We will use the notation of [\ref{Mil:IC}] and denote
the $n-th$ iterate of a map $ f $ by $ f^{\circ n} $. Let $ {\bf
IR}^{2} $ be the set of all intervals in $ {\bf R}^{2} $. Let $
If: {\bf IR}^{2} \rightarrow {\bf IR} $ be an interval extension
of $ f $. Let $ If^{\circ n} $ be the $ n-th $ iterate of the map
$ If $ . We can consider the evolution of $ If $  as (i) evolution
of two correlated real dynamical systems $[{\underline If}^{\circ
n},{\overline If}^{\circ n}],  $ or as (ii) evolution of an
interval dynamical system on interval data. The principal
considerations concerns case (i). In the case for iterated map $$
[{\underline x_{k+1}}, \; {\overline x_{k+1}}] =
 [{\underline f}({\underline x_{k}},{\overline x_{k}},
 {\underline u},{\overline u},{\underline v},{\overline v}), \;
 {\overline f}({\underline x_{k}},{\overline x_{k}},
 {\underline u},{\overline u},{\underline v},{\overline v})) ] $$
in $ {u,v} $ plane we will describe the evolution of $ If $ for
various subregions and points of $2-$dimensional region
 $$ D_{p}: \; \infty > u > 1, \; 1 \leq v \leq
(2^u - 1)^{\frac{1}{u}}, $$
where $ v $ is some real parameter.
Let us give formulas for interval extension of Minkowski's moduli space:
here the bounds $ \; {\underline \tau},
 \; {\overline \tau}, \; $ are obtained with the help of the iteration
process:

$$ {\underline t}_{i + 1} = (1 + {\underline t}_{i}^{\overline
p})^ {\frac{1}{\overline p}}(((1 - (1 + {\underline
t}_{i}^{\overline p})^ {-\frac{1}{\overline p}} - (1 + {\overline
\sigma}^{\underline p})^ {-\frac{1}{\underline p}})^{\underline
p})^{\frac{1}{\underline p}} - {\overline \sigma}(1 + {\overline
\sigma}^{\underline p})^{-\frac{1} {\underline p}})~\eqno (1) $$

 $${\overline t}_{i + 1} = (1 + {\overline t}_{i}^{\underline p})^
{\frac{1}{\underline p}}(((1 - (1 + {\overline t}_{i}^{\underline
p})^ {-\frac{1}{\underline p}} - (1 + {\underline
\sigma}^{\overline p})^ {-\frac{1}{\overline p}})^{\overline
p})^{\frac{1}{\overline p}} - {\underline \sigma}(1 + {\underline
\sigma}^{\overline p})^{-\frac{1} {\overline p}})~\eqno (2) $$

    $$         \; i = 0,1,\cdots  $$
As interval computation is the enclosure method, we have to put:

$$ [{\underline \tau}, \; {\overline \tau}] =
 [{\underline t}_{N}, \; {\overline t}_{N}] \bigcap
 [{\underline \tau}_{0}, \; {\overline \tau}_{0}] \; .$$
$ N $ is
computed on the last step of the iteration. \\
For initial values we may take $: \; [{\underline t}_{0}, \;
{\overline t}_{0}] =
[{\underline \tau}_{0}, \; {\overline \tau}_{0}] = [0,\; 0.36]. $
%\begin{expl}
%  Computation on intervals (trajectories in interval spaces).
%\end{expl}
\begin{remk}
Let $ {\bf X} = [{\underline p}, {\overline p}; \; {\underline
\sigma}, {\overline \sigma}]$ be the interval, where the interval
iteration process~$(1) - (2)$ is computed. Let $p =
\frac{{\underline p} + {\overline p}}{2}, \sigma =
\frac{{\underline \sigma} + {\overline \sigma}}{2}$. Computations
show that for the convergence of the interval iterative
processes~$(1) - (2)$ it is sufficient that the noninterval
inequality $f_{\tau}^{'} < 1$ is satisfied in the point
$[p,\sigma].$
\end{remk}

\section{Dynamical systems from critical lattices}

\subsection{Algebraic dynamical systems}
  Let $\Lambda $ be a lattice and $R$ its ring of multipliers.
  So for $\lambda \in R, \; \lambda\Lambda \subseteq \Lambda$ and $\omega \in \Lambda$
  an algebraic ${\bf Z}-$action
  $\alpha: n \mapsto \alpha_{n}$ on $\Lambda $ is defined by
  $$
   \alpha_{n}(\omega) = {\lambda}^{n}\omega.
  $$
\subsection{Elliptic curves from critical lattices}
Let $\Lambda $ be a critical lattice of the domain $ D_{p}$.
 After complexification of ${\bf R}^2$ the lattice $\Lambda$ takes
 form  $\Lambda = n\omega_1 + m\omega_2, \; \omega_1, \omega_2 \in {\bf C},
  n, m \in {\bf Z}, \frac{\omega_1}{\omega_2}$ is not a real
  number.
\begin{expl}
After complexification the basis of the critical lattice
$\Lambda_{2}^{(0)}$ has the form  $\omega_1 = 1, \;
  \omega_2 = 1/2 + \frac{\sqrt 3}{2}i.$ \\
  Respectively for $2 \le p \le p_{0}, $ the lattice $\Lambda_{p}^{(0)}$ has the
  basis  $\omega_1 = 1, \;
  \omega_2 = 1/2 + \frac{\sigma_{p}}{2}i.$
\end{expl}
 Let $\Lambda $ be as above. For $\alpha \in \Lambda$ construct
 invariants of $\Lambda :$
 $$ \sum^{'}_{\alpha \in \Lambda} \frac{1}{\alpha^n} =
 \sum_{\alpha \ne 0, \; \alpha \in \Lambda} \frac{1}{\alpha^n} .$$
 Let $$ c_{n} = \sum^{'}_{\alpha \in \Lambda}
 \frac{1}{\alpha^{2n}}.$$ There is well known
\begin{lem}
 If  $t \in  {\bf R}, t > 2,$ then the series
 $ \sum^{'}_{\alpha \in \Lambda} \frac{1}{\alpha^{t}},
 $
 converges absolutely.
\end{lem}

  The Weierstrass elliptic function is defined as expression
  $$ \frac{1}{z^2} + \sum^{'} (\frac{1}{(z + \alpha)^{2}} -
  \frac{1}{\alpha^{2}}).
  $$

% \begin{expl}
 % Weierstrass elliptic functions of critical lattices for $p = 2.$
%\end{expl}

Each Weierstrass elliptic function defines the field ${\mathcal
K}_{\Lambda}$ of elliptic functions. Its model is the elliptic
curve in the Weierstrass normal form: $$ y^2 = 4x^3 - 60c_{2}x -
140c_{3}. $$

\subsection{Dynamical systems from the elliptic curves}
 The {\it Julia set} of a rational function $f \in {\bf C}(z)$ is
 the closure of the union of its repelling cycles [\ref{Mil:IC}].
 Let $f(z)$ be the rational function that is constructed under division of points of
 an elliptic curve over ${\bf C} [\ref{S:EC}].$
 The Julia set of the rational functions is by the result of S. Latt{\'e}s
 [\ref{L:IS}] the whole  sphere  $\overline {\bf C}.$

%\vskip 15 pt

\vskip 10 pt

\centerline{\Large REFERENCES} \vskip 7 pt

\begin{enumerate}
\small

\item
\label{V:CW}
  G.F. Vorono{\"i} (1952),
   {\it Collected works in 3 volumes},
   Kiev: Acad.of Sci. Ukr.SSR.

\item
\label{C:GN} J. Cassels (1971) {\it An Introduction to the
Geometry of Numbers}, Berlin: Springer-Verlag.

\item
\label{Mi:DA} H. Minkowski (1907), {\it Diophantische
Approximationen}, Leipzig: Teubner.

\item
\label{GGM:PM} N. Glazunov, A. Golovanov, A. Malyshev (1986), {\it
Proof of Minkowski
     hypothesis about critical determinant of
$ |x|^p + |y|^p <1 $      domain},
 Research in the number theory.9. Notes of scientific seminars
of LOMI. {\bf 151} Leningrad: Nauka.  40--53.

\item
\label{D:AN} C. Deninger (1999), {\it Some analogies between
number theory and dynamical systems on foliated spaces},  Doc.
Math. J. DMV Extra volume ICM1998. 163--186.

\item
\label{DEMW:DS} P. D'Ambros, G. Everest, R. Miles, T. Ward (2000),
 {\it Dynamical systems arising from elliptic curves},
   Colloq. Math.,  {\bf 84-85}, Pt. 1,  95--107

\item
\label{M:LP} L.J. Mordell (1941), {\it Lattice points in the
region} $ \mid Ax^4 \mid + \mid By^4 \mid \geq 1,$   J. London
Math. Soc. {\bf 16} , 152--156.

\item
\label{D:NC} C. Davis (1948), {\it Note on a conjecture by
Minkowski},  Journ. of the London Math. Soc.,  {\bf 23}, (3),
172--175.

\item
\label{Co:MC} H. Cohn (1950), {\it Minkowski's conjectures on
critical lattices in the metric \\ $ \{\mid \xi \mid^p + \mid \eta
\mid^p \}^\frac{1}{p},$}  Annals of Math.,  {\bf 51}, (2),
734--738.

\item
\label{W:MC} G. Watson (1953), {\it Minkowski's conjecture on the
critical lattices of the region} $ |x|^p + |y|^p \leq 1 \;$, (I),
(II),  Journ. of the London Math. Soc.,  {\bf 28}, (3, 4),
305--309, 402--410.

\item
\label{Ma:AC} A. Malyshev (1977), {\it Application of computers to
the proof of a conjecture of Minkowski's from geometry of numbers.
1},  Zap. Nauchn. Semin. LOMI,  {\bf 71},  163--180.

\item
\label{AH:IC} G. Alefeld, J. Herzberger (1983), {\it Introduction
to Interval Computations}, NY: Academic Press, .

\item
\label{G:IA} N.M. Glazunov (1997),   {\it On Some
Interval-Algebraic Methods for Verifications of Dynamical
Systems}, Cybernetics and Computer Technologies,  {\bf 109},
15--23.

\item
\label{Mil:IC}
 J. Milnor (1992), {\it Remarks on Iterated Cubic Map},  Experimental
Mathematics, {\bf 1} , No. 1, 5--24.

\item
\label{G:RM} N. Glazunov (2001), {\it Remarks to the Minkowski's
conjecture on critical determinant of the region}, $ |x|^p + |y|^p
<1, p > 1,$  Theses of the reports
           to the 4-th Int. Conf. on Geometry and Topology,
           Cherkasi:   ChITI,  21--23.

\item
\label{S:EC} J. Silverman (1986), {\it The Arithmetic of Elliptic
Curves}, NY: Springer.

\item
\label{L:IS} S. Latt{\'e}s (1918), {\it Sur l't{\'e}ration des
substitutions rationelles et les fonctions de Poincar{\'e}}, C. R.
Acad. Sci. Paris,  {\bf 166}, (3), 26--28.

\end{enumerate}

\vskip 20 pt

\noindent {\it Glushkov Institute of Cybernetics NAS }

\vskip 2 pt

\noindent {\it e-mail: {\tt glanm@yahoo.com}}

\end{document}